\documentclass{article}

\newtheorem{example}{Example}
\newtheorem{lemma}{Lemma}
\newtheorem{theorem}{Theorem}
\newtheorem{remark}{Remark}

\usepackage{amsfonts}
\usepackage{amsmath}
\usepackage{amssymb}
\usepackage{graphicx}

\numberwithin{equation}{section}

\begin{document}

\title{The Finite-Time Turnpike Phenomenon for Optimal Control Problems: Stabilization by
Non-Smooth Tracking Terms}

\author{Martin Gugat, Michael Schuster, Enrique Zuazua}
%
%
\date{
Lehrstuhl f\"ur Angewandte Analysis (Alexander-von-Humboldt Professur), 
 Department Mathematik,
 Friedrich-Alexander Universit\"at Erlangen-N\"urnberg (FAU),
  Cauerstr. 11, 91058 Erlangen, Germany.
  This is a preprint of the corresponding chapter
  in the book 
  "Stabilization of Distributed Parameter Systems: Design Methods and Applications" within the SEMA SIMAI Springer Series edited by Alexander Zuyev 
  reproduced with the permission of
  the Publisher as it appears on the copyright page of the book.
  }

%


\maketitle




\abstract{
In this paper,  problems of optimal  control
are considered where
in the objective function,
in addition to the control cost there is
a
tracking term
that
measures the distance to a
desired stationary state.
The tracking term is given by
some norm and therefore it is in
general not differentiable.
In the optimal control problem,
the initial
state
is prescribed.
We assume that the system is
either exactly controllable in
the classical sense or nodal profile controllable.
%
We show that both for
systems that are governed by ordinary differential equations
and for infinite-dimensional systems,
for example for boundary control systems governed by
the wave equation,
under certain assumptions
the optimal system state
is steered exactly to
the desired state after finite time.
}


\section{Introduction}
Since the turnpike  phenomenon  has been studied by P. A. Samuelson in mathematical economics in 1949
(see \cite{samuelson}),
it has been analyzed in various contexts, see for example
 \cite{zsalavski2000},
 \cite{zsalavskiband2006}
and \cite{ANEXPONENTIAL}.
For  optimal control problems with partial differential equations
it has been studied in
\cite{porretta}
and
\cite{trelatzhangzuazua}
where distributed control is considered for
 linear--quadratic optimal control problems.
 Problems of optimal boundary control
 are studied in \cite{guzu}, \cite{guha} and \cite{gupafa}.
In \cite{trelatzhang}, both
integral- and measure--turnpike properties
are considered.
The turnpike phenomenon for
linear quadratic optimal control problems
with time-discrete systems is studied in
\cite{gugli}.
In \cite{schaller},
 linear quadratic optimal control problems governed by general evolution
equations are considered
and exponential sensitivity and turnpike analysisis studied.
An overview on the turnpike phenomenon is
given in the monograph \cite{zas2019}.

In this paper, we consider
integral turnpike  properties
for problems where
the system is exactly controllable
and in the objective function,
an $L^1$-norm or $L^2$-norm tracking term appears.
We show that the resulting optimal
controls have a finite-time  turnpike structure,
that is
the optimal state reaches the
static desired state (that we also refer to as the {\em turnpike}
and that does not depend on time)
exactly  in finite time.
%

These turnpike result are also useful for numerical computations
since they show that
for sufficiently large time horizons $T$,
sufficiently accurate
approximations of the optimal  state/control pairs
should  also  be
identical to the desired state
with the corresponding constant control
most of the time.

The finite-time (or exact)  turnpike property
for continuous-time systems
has already   been discussed in
   \cite{faulwasserIFAC}
   as an assumption
   in the context of
   nonlinear model predictive control
   for a finite-dimensional system
   that is governed by an ordinary differential
   equation.
   Here the aim is to prove convergence in
   model predictive control.
   As an application, a problem of optimal fish harvesting control is studied.

This paper has the following structure.
In order to illustrate the situation,  first we
consider optimal control problems that are governed by
ordinary differential equations.
In  these problems the $L^1$--norm appears in
the tracking term in the objective function.
We show that
if the weight  of the tracking term  (i.e. the penalty parameter)
is sufficiently large,  the optimal
states and controls
have a finite-time   turnpike structure.

In the next section,
we present a finite-time turnpike result
for optimal control problems
with an abstract infinite dimensional system.
First we consider the case where the system is
exactly controllable.
We consider an optimal control problem
where the tracking term is given by
a certain maximum norm.
We show that if the weight of the tracking
term is sufficiently large,
the solution has a finite-time turnpike structure.

Then we consider the case where the system is
nodal profile exactly controllable.
We consider an optimal control problem
where the tracking term for the nodal profiles is given by
an  $L^2$--norm.
We show that if the weight of the tracking
term is sufficiently large,
the solution has a finite-time turnpike structure
for the nodal profiles.

Finally we return to the case where the system is
exactly controllable.
We consider an optimal control problem
where the tracking term is given by
a weighted $L^1$--norm that has a singularity at $t_0>0$.
We show also in this case,
the solution has a finite-time turnpike structure.

In Section   \ref{example}, examples are presented,
where the results from the previous section are applicable.
Section \ref{conclusion} contains conclusion.


\section{Optimal control problems with ordinary differential equation}
We start with optimal control problems
with systems that are governed by ordinary differential equations.
We show that for such systems, $L^1$-tracking terms
in the objective function can lead to
finite-time turnpike structures.

%

\begin{example}
\label{bsp1}
We start with a system similar to the motivating example in \cite{guzuexact}
that is governed by an ordinary differential equation.
Let $\gamma >0 $ be given.
For $T>0 $ sufficiently large (this will be specified later) we consider the problem
\[
{\bf (OC)}_T
\left\{
\begin{array}{l}
\min\limits_{u\in L^2(0,T)} \;
\int\limits_0^T \,
\tfrac{1}{2} |u(t)|^2 + |u(t)| + \gamma \,  |y(t)|\, dt
\; \,{\sf subject } \;{\sf to}\;
\\
\\
y(0)= -1, \;y'(t)= y(t) + \exp(t) \, u(t).
\end{array}
\right.
\]
The corresponding  optimal control problem
where  the initial condition does not appear is
\[
{\bf (OC)}^{(\sigma)}
\left\{
\begin{array}{l}
\min\limits_{u\in L^2(0,T)} \;
\int\limits_0^T \,
\tfrac{1}{2} |u(t)|^2 + |u(t)| + \gamma \,  |y(t)|\, dt
\; \,{\sf subject } \;{\sf to}\;
\\
\\
y'(t)= y(t) + \exp(t) \, u(t).
\end{array}
\right.
\]
The solution of  ${\bf (OC)}^{(\sigma)}$
(that we call the turnpike)
 is zero, that is
$y^{(\sigma)}=0$
and
$u^{(\sigma)}=0$.
The results about the solution of
${\bf (OC)}_T$ are summarized in the following lemma.
\begin{lemma}
\label{bsplemma1111}
For $\gamma >0 $, define $t_0>0$ as the minimal value where
\[
(t_0 - 1) \, \exp(t_0) = \frac{1}{\gamma } - 1. 
\]
Assume that $T>t_0$ and (even)
\begin{equation}
\label{11112019}
 \gamma \, {\rm e}^T \geq  1 + \gamma  \, {\rm e}^{t_0}.
 \end{equation}
Define
\begin{equation}
\label{hatudefinition}
\hat u(t) = \gamma  ({\rm e}^{t_0} -  {\rm e}^t)\geq 0  \;\; \mbox{\rm  for }\;\; t\in (0,\, t_0],\;
\hat u(t)=0
 \;\; \mbox{\rm   for }\;\;t > t_0.
 \end{equation}
 Then for the state $\hat y$ generated by $\hat u$ for $t\geq t_0$ we have $\hat y(t) = 0$.
 Moreover, for all $ t \in (0, T)$ we have $\hat y(t)\leq 0$.

 The control   $\hat u$ as defined in (\ref{hatudefinition})  is the
 unique solution of ${\bf (OC)}_T$.

\end{lemma}
{\bf Proof.}
Let a control $u\in L^2(0, \, T)$ be given.
Then for the corresponding state $y$
we have
\begin{equation}
\label{zustandy}
y(t) = {\rm e}^t \left[ - 1 + \int_0^t u(\tau)\, d\tau \right]
.
\end{equation}
Note that for the optimal control we have
$y(t)\leq 0$.
(If $y(t_0)=0$, we can continue with the zero control.)
Moreover, we have $u(t) \geq 0$.
(Otherwise, instead of decreasing the state it is also better
to switch off the control).
Hence it suffices to consider the feasible controls $u(t) \geq 0$
that satisfy the moment inequality
\begin{equation}
\label{momentengleichung}
\int_0^T u(\tau)\, d\tau \leq 1.
\end{equation}
Due to the definition of $t_0$ and (\ref{hatudefinition})  we have
\begin{equation}
\label{integralgleichung}
\int_0^T \hat  u(\tau)\, d\tau  = 1.
\end{equation}

Then for $t\in (0, t_0)$ we have
\[
\hat y(t) = {\rm e}^t \left[ - 1 + \int_0^t u(\tau)\, d\tau \right]
=
 \gamma \, t \,  {\rm e}^{t+t_0} -   \gamma \, {\rm e}^{2t} + (\gamma - 1){\rm e}^{t} \leq 0
\]
and for $t\geq t_0$ we have $\hat y(t) = 0$.

With the representation (\ref{zustandy}),
for all feasible controls $u\geq 0$ where $y\leq 0$
integration by parts yields
\[J_{(0, T)}(u,\, y)
=
\int\limits_0^T \,  \left[  \tfrac{1}{2} |u(  t)|^2 + u(t)  - \gamma \,  y(t)\right] \, dt
\]
\[
=
\int\limits_0^T \,
   \left\{
   \tfrac{1}{2} |u(  t)|^2 + u(t)  + \gamma \, {\rm e}^t \left[ 1 -  \int_0^t u(\tau)\, d\tau \right]
   \right\}
   \, dt
\]
\[
=
\int\limits_0^T \,u(t) \, dt
 + \int\limits_0^T \,   \tfrac{1}{2} |u(  t)|^2
 +
\gamma \, {\rm e}^t \left[ 1 -  \int_0^t u(\tau)\, d\tau \right]|_{t=0}^T
+
\int\limits_0^T \, \gamma \,{\rm e}^t \,  u(t) \, dt
\]
\[
=
\int\limits_0^T \,u(t) \, dt
- \gamma +
\gamma \, {\rm e}^T \left[ 1 -  \int_0^T u(\tau)\, d\tau \right]
+
\int\limits_0^T \,  \left[  \tfrac{1}{2} |u(  t)|^2
+
 \gamma \,{\rm e}^t \,  u(t) \right]\, dt
\]
\[
=
( \gamma \, {\rm e}^T - 1) \, \left[ 1 -  \int_0^T u(\tau)\, d\tau \right]
+
\int\limits_0^T \,  \left[ \tfrac{1}{2} |u(  t)|^2
+
\gamma \, {\rm e}^t \,  u(t) \right] \, dt
+
1  - \gamma.
\]

If $T$ is sufficiently large in
the sense that  (\ref{11112019}) holds,  due to the $L^1$-norm that appears in the objective function,
the solution has an exact turnpike structure where
the system is steered to zero in  the  finite time
$t_0$ that is independent of $T$ and remains
there for $t\in (t_0,\, T)$.
This can be seen as follows.
Let $u(t) = \hat u(t) + \delta(t)$
with $\hat u$ as defined in (\ref{hatudefinition}) and
$\int_0^T \delta(\tau) \, d\tau \leq  1 -  \int_0^T \hat  u(\tau)\, d\tau  = 0$
where the last equation follows from
(\ref{integralgleichung})
and $\delta(t) \geq 0$ for $t\geq t_0$.
Due to (\ref{11112019}) we have

\[J_{(0, T)}(u,\, y)
=
( \gamma \, {\rm e}^T - 1)  \left[  -  \int_0^T \delta(\tau)\, d\tau \right]
\]
\[
+
\int\limits_0^T \, \left[   \tfrac{1}{2} |\hat u(  t) + \delta(t)|^2
+
\gamma \,  {\rm e}^t \, (\hat u(  t) + \delta(t))\right] \, dt+
1  - \gamma
 \]
 \[
 \geq
 ( \gamma \, {\rm e}^T - 1)  \left[  -  \int_0^T \delta(\tau)\, d\tau \right]
 +
 \int\limits_0^T \,   \tfrac{1}{2} \hat u(  t)^2  +
 \gamma \, {\rm e}^t \, \hat u(  t) \,dt
 +
 \int\limits_0^T \, \left( \hat u(  t) +\gamma\, {\rm e}^t \right) \,
 \delta(t) \, dt
 +
1  - \gamma
\]
\[
=
( \gamma \, {\rm e}^T - 1)  \left[  -  \int_0^T \delta(\tau)\, d\tau \right]
+
J_{(0, T)}(\hat u,\, \hat y)
 + \int\limits_0^{t_0}  \,
\gamma {\rm e}^{t_0}  \,  \delta(t) \, dt
+
\int\limits_{t_0}^T  \gamma \,  {\rm e}^t  \, \delta(t) \, dt
\]
\[
\geq
( \gamma \, {\rm e}^T - 1)  \left[  -  \int_0^T \delta(\tau)\, d\tau \right]
+
J_{(0, T)}(\hat u,\, \hat y)
 +
\int\limits_0^{t_0}  \gamma \,  {\rm e}^{t_0}  \,  \delta(t) \, dt
+
\int\limits_{t_0}^T  \gamma \,   {\rm e}^{t_0}  \, \delta(t) \, dt
\]
\[
=
J_{(0, T)}(\hat u,\, \hat y)
+
( \gamma \, {\rm e}^T - \gamma \,  {\rm e}^{t_0} - 1)
\left[  -  \int_0^T \delta(\tau)\, d\tau \right]
.
\]
Since $\gamma \, {\rm e}^T - \gamma \,  {\rm e}^{t_0} - 1 \geq 0$,
this implies that  $\hat u$ as defined in (\ref{hatudefinition})  is the optimal control.
Thus we have proved Lemma \ref{bsplemma1111}. $\Box$

Consider the value $t_0$ as a function of
$\gamma$, $t_0=t_0(\gamma)$.
Then we have $t_0( 1 ) = 1$
and
\[
\lim_{\gamma \rightarrow \infty} t_0(\gamma) = 0.
\]
In Example \ref{beispielnumerik},
we present numerical approximations
for the optimal states and controls
for three values of $\gamma$.
%
%
%
\end{example}


\subsection{A  more general result for scalar ordinary differential equations}
Now we consider an optimal control problem
with  the same objective function and
a  more general  ordinary differential equation.
In this problem, we also prescribe a terminal condition.
At the end of the section
we will  present sufficient conditions
that  imply that
if the penalty parameter $\gamma$ is sufficiently large,
the terminal state is reached
before  the final time.

%
Let continuous
functions $f$, $g$ from $[0, \infty)$ to
the real numbers be given.
Assume that for all $t\geq 0$ we have $f(t)>0$, and
$g(t)>0$.
Let $\gamma \geq 1$ and $\alpha < 0$ be given.
For  a  finite time horizon $T > 0$ we consider the problem
\[
{\bf (OC)}_T
\left\{
\begin{array}{l}
\min\limits_{u(t)\in L^2(0, \, T), y(t)\in AC(0, \, T)} \;
\int\limits_0^T \,
\tfrac{1}{2} |u(t)|^2 + |u(t)| + \gamma \,  |y(t)|\, dt
\; \,{\sf subject } \;{\sf to}\;
\\
y(0)= \alpha, \;y'(t)= f(t) \, y(t) + g(t) \, u(t)
\\
y(T)= 0.
\end{array}
\right.
\]
Here again the solution
of the corresponding optimal control problem
without the initial and the terminal conditions
(the turnpike) is zero, that is
$y^{(\sigma)}=0$
and
$u^{(\sigma)}=0$.
Note that the turnpike is compatible
with the terminal constraint $y(T)=0$.
In the following theorem we present the
optimal control for
${\bf (OC)}_T$,
which  has a similar structure as
in the previous example.
\begin{theorem}
Define
\[F(t) = \exp(\int_0^t f(s) \, ds),
\;
H(t) = \int_0^t F(\tau) \, d\tau.
\]
We have
\begin{equation}
\label{yrepresentation24092019}
y(t) = F(t) \, \left[
\alpha + \int_0^t \frac{g(\tau)}{F (\tau)} \, u(\tau) \, d \tau\right].
\end{equation}
Define
\begin{equation}
\label{optimalcontroldefinition}
\hat u(t) =
\max\left\{0,\;
 \left[ -1 -\gamma \, \frac{g(t) \, H(t)}{F(t)} + \lambda \frac{g(t)}{F(t)} \right]
 \, \right\}
\end{equation}
where the number $\lambda > 0$ is chosen such that
\begin{equation}
\label{lambdawahlbedingung}
\int_0^T \hat u(\tau) \, \frac{g(\tau)}{F(\tau)} \, d\tau = - \alpha.
\end{equation}
Then the unique optimal control that solves
${\bf (OC)}_T$ is equal to $\hat u(t)$.
\end{theorem}

{\bf Proof.}
Since $y(0) = \alpha \leq 0$,
for the optimal state we have $y(t) \leq 0$
for all $t\geq 0$.
(Since otherwise, instead of increasing the
state above zero
it is better to switch off the control.)
Moreover, for the optimal control we have
$u(t) \geq 0$.
(Since otherwise, instead of decreasing
the state it is also better to switch off the control).
Hence it suffices to consider the feasible controls $u(t) \geq 0$
that satisfy the moment inequality
\begin{equation}
\label{momentengleichungallgemeinesbeispiel}
\int_0^T \frac{g(\tau)}{F (\tau)} \, u(\tau)\, d\tau \leq - \alpha.
\end{equation}
Due to the choice of $\lambda$,
for the state $\hat y$ generated by $\hat u$, we have
$\hat y(T)=0$.
For $t\in [0, \, T]$, consider
\begin{equation}
\label{Bdefinition}
B(t) = \int_0^t  \frac{g(\tau)}{F(\tau)} \, \hat u (\tau)   \, d\tau.
\end{equation}
Then $B(0)=0$ and $B$  is increasing.
Hence also the function
$[\alpha + B(t)]$ is increasing.
We have  $B(0)+ \alpha  < 0$ and $B(T)  + \alpha =0$.
Thus there exists a  unique point
\[t_0 = \min\{t\in [0, \, T]: \alpha + B(t)=0\}.\]
and we have $t_0 \in (0,\, T]$.
We have $B(t_0)= B(T)$ and $B$ is increasing.
This implies
that for all $t\in [t_0, \, T]$,
we have
$B(t) = - \alpha$.
On account of the definition of $B$ as
an integral, this is only possible
if for  all $t\in [t_0, \, T]$, we have $\hat u(t)=0$.
This implies that for  all $t\in [t_0, \, T]$ we have
\begin{equation}
\label{ungleichunggroessertnull}
-1 -\gamma \, \frac{g(t) \, H(t)}{F(t)} + \lambda \frac{g(t)}{F(t)} \leq 0.
\end{equation}
By (\ref{yrepresentation24092019}) we have
\[t_0 = \min\{t\in [0, \, T]: \hat y(t)=0\}.\]

Since $\hat y(t) = F(t) \left[ \alpha + B(t)\right]$,
for $t < t_0$ we have $\hat y(t) < 0$.
Since for $t \geq t_0$, we have $\hat u(t) = 0$,
this implies that  $\hat y(t) =  0$ for all $t\geq t_0$.

Since $\hat u \geq 0$ and $\hat y \leq 0$,
for the objective function we have
\begin{eqnarray*}
J(\hat u) & = & \int_0^{ T  }
\frac{1}{2} \,
(\hat u(t))^2 + \hat u(t) -
\gamma \,  F(t) \left[
\alpha + \int_0^t \hat u (\tau)
\,
\frac{g(\tau)}{F(\tau)}\, d\tau \right] \, dt.
\end{eqnarray*}
Integration by parts yields (since $B(T)= - \alpha$)
\begin{eqnarray*}
J(\hat u) & = &
\int_0^{T } \frac{1}{2} \,
(\hat u(t))^2 + \hat u(t) \, dt
- \gamma \, H(s) \left[
\alpha +  \int_0^s \hat u (\tau)
\,
\frac{g(\tau)}{F(\tau)}\, d\tau \right]|_{s=0}^{T}
\\
&
+
&
\gamma \,\int_0^{T} H(t) \,  \hat u (t)
\,
\frac{g(t)}{F(t)}\,dt
\\
& = &
\int_0^{ T } \frac{1}{2} \,
(\hat u(t))^2 + \hat u(t) \, dt
 - \gamma \,H(t_0) \,  (\alpha + B(T) )
 \\
 &  + &\gamma \,
\int_0^{T} \hat u(\tau) \, \frac{ H(\tau) \, g(\tau)}{ F(\tau)} \, d\tau
\\
& = &
\int_0^{T } \frac{1}{2} \,
(\hat u(t))^2 + \hat u(t)
\left[ 1 + \gamma \,\frac{ H(t) \, g(t)}{ F(t)} \right] \, dt.
\end{eqnarray*}
Let $\delta \in L^2(0,\, T)$ be given.
We use $\delta$ as a perturbation of
the control.
To make sure that the terminal  condition remains valid,
we assume that
\begin{equation}
\label{momentlambda}
\int_0^{T} \delta(\tau) \, \frac{g(\tau)}{F(\tau)} \, d\tau = 0.
\end{equation}
Since the optimal control
must increase the values of
the corresponding trajectory  to zero,
it can only have positive values.
Therefore we  assume that for $t\in [0, \, T]$ we have
$\hat u(t) + \delta(t)\geq 0$. Thus
 for $t\in [0, \, t_0]$ we have
${\rm sign}(\hat u(t) + \delta(t) ) =1$
and for $t\geq t_0$, we have $\delta(t) \geq 0$.
Then we have
\begin{eqnarray*}
J(\hat u + \delta) & = &
\int_0^{T } \frac{1}{2} \,
(\hat u(t) + \delta(t) )^2 + (\hat u(t)+ \delta(t))
\left[ {\rm sign}(\hat u(t) + \delta(t) ) + \gamma \, \frac{ H(t) \, g(t)}{ F(t)} \right] \, dt
\\
& = &
J(\hat u)
+ \int_0^{T }  \frac{1}{2}  \delta(t)^2  \, dt
+
 \int_0^{ t_0 }
 \delta(t) \, \left[\hat u +  1 + \gamma \,\frac{ H(t) \, g(t)}{ F(t)} \right] \, dt
 \\
 & + &
  \int_{ t_0 }^T
 \delta(t) \, \left[   {\rm sign}(\delta(t) )  + \gamma \,\frac{ H(t) \, g(t)}{ F(t)} \right] \, dt
\\
& = &
J(\hat u)
+ \int_0^{T }  \frac{1}{2}  \delta(t)^2 \, dt
+  \int_0^{ t_0 }
\delta(t) \, \lambda \, \frac{  g(t)}{F(t)} \, dt
\\
&
+
&
\int_{ t_0 }^T
 \delta(t) \, \left[   {\rm sign}(\delta(t) )  + \gamma \,\frac{ H(t) \, g(t)}{ F(t)} \right] \, dt
\\
& = &
J(\hat u)
+ \int_0^{T }  \frac{1}{2}  \delta(t)^2 \, dt
+  \int_0^{ T }
\delta(t) \, \lambda \, \frac{  g(t)}{F(t)} \, dt
\\
&
+
&
\int_{ t_0 }^T
 \delta(t) \, \left[
 - \lambda \, \frac{  g(t)}{F(t)} +
   {\rm sign}(\delta(t) )  + \gamma \,\frac{ H(t) \, g(t)}{ F(t)} \right] \, dt
\\
& = &
J(\hat u)
+ \int_0^{T}  \frac{1}{2}  \delta(t)^2\, dt
\\
&
+
&
\int_{ t_0 }^T
 \delta(t) \, \left[
   1  + \gamma \,\frac{ H(t) \, g(t)}{ F(t)}
 - \lambda \, \frac{  g(t)}{F(t)}
 \right] \, dt
 \\
 &\geq &
 J(\hat u)
\end{eqnarray*}
where the last step follows with (\ref{ungleichunggroessertnull}).
Thus $\hat u$ is the minimizer of $J$
among  all controls that generate states with $y(T)=0$.
This shows the assertion. $\Box$

The question remains:
Do we have $t_0 < T$ if $\gamma$ is sufficiently large?

Let $t_1\in (0, \, T)$ be given
such that
\begin{equation}
\label{t1condition}
-\alpha  -  \frac{F(t_1)}{g(t_1)} \, \int_0^{t_1 } \frac{g^2}{F^2} \, dt
+
\int_0^{t_1}
\frac{ g} {F} \,dt
>
0
.
\end{equation}

Note that $H$ is strictly  increasing,
hence we have the inequality
 \[\int_0^{t_1 } \left(H(t_1)  - H(t) \right) \frac{g^2}{F^2}  \, dt >0.\]
Define  the number
\begin{equation}
\label{gammadefinition}
\gamma(t_1) =
\frac{ -\alpha  -  \frac{F(t_1)}{g(t_1)} \, \int_0^{t_1 } \frac{g^2}{F^2} \, dt
+
\int_0^{t_1}
\frac{ g} {F} \,dt }
{ \int_0^{t_1 } \left(H(t_1)  - H(t) \right) \frac{g^2}{F^2}  \, dt }.
\end{equation}
Then we have $\gamma(t_1)>0$.
Define the number
\begin{equation}
\label{lambdadefinition}
\lambda_1 =
\frac{-\alpha
+
\int_0^{t_1}
\frac{ g} {F} \,dt  + \gamma(t_1) \int_{0}^{t_1} H  \frac{g^2}{F^2}  \, dt}{
\int_0^{t_1 }
\frac{g^2}{F^2} \, dt
}
.
\end{equation}
The  definition of $\lambda_1$
implies  the equation
\begin{equation}
\label{momentequationt0}
\lambda_1 \,
\int_0^{t_1 }
\frac{g^2}{F^2} \, dt =
-\alpha
+
\int_0^{t_1}
\frac{ g} {F} \,dt  + \gamma(t_1) \int_{0}^{t_1} H  \frac{g^2}{F^2}  \, dt.
\end{equation}
Moreover, due to the definition of $\gamma(t_1) $ we  have
\[
1
+
\frac{g(t_1) \, H(t_1)}{F(t_1)} \,
\gamma(t_1) =
\frac{
\int_{0}^{t_1} H  \frac{g^2}{F^2}  \, dt
+\frac{g(t_1) \, H(t_1)}{F(t_1)}
\left[\alpha
-
\int_0^{t_1}
\frac{ g} {F} \,dt
\right]
}
{\int_{0}^{t_1} H  \frac{g^2}{F^2}  \, dt
- H(t_1) \, \int_0^{t_1 }
\frac{g^2}{F^2} \, dt }
 .
 \]
In addition, the definition of $\lambda_1$ and of  $\gamma(t_1)$ implies
\[
\lambda_1 \, \frac{g(t_1)}{F(t_1)}  =
\frac{
-
\frac{g(t_1)}{F(t_1)}
\left[
\alpha
-
\int_0^{t_1}
\frac{ g} {F} \,dt
\right]
  +
  \frac{g(t_1)}{F(t_1)}
  \gamma(t_1) \int_{0}^{t_1} H  \frac{g^2}{F^2}  \, dt}{
\int_0^{t_1 }
\frac{g^2}{F^2} \, dt
}
\]
\[
=
\frac{
-
\frac{g(t_1)}{F(t_1)}
\left[
\alpha
-
\int_0^{t_1}
\frac{ g} {F} \,dt
\right]
  +
 \left[
 \frac{
 \int_0^{t_1 }
\frac{g^2}{F^2} \, dt
+\frac{g(t_1) }{F(t_1)}
\left[\alpha
-
\int_0^{t_1}
\frac{ g} {F} \,dt
\right]
}
{\int_{0}^{t_1} H  \frac{g^2}{F^2}  \, dt
- H(t_1) \, \int_0^{t_1 }
\frac{g^2}{F^2} \, dt }
  \right]
  \int_{0}^{t_1} H  \frac{g^2}{F^2}  \, dt}{
\int_0^{t_1 }
\frac{g^2}{F^2} \, dt
}
\]
\[
=
\frac{1}{
\int_0^{t_1 }
\frac{g^2}{F^2} \, dt
}
\left[
\frac{
\int_0^{t_1 }
H\,
\frac{g^2}{F^2} \, dt
}
{\int_{0}^{t_1} H  \frac{g^2}{F^2}  \, dt
- H(t_1) \, \int_0^{t_1 }
\frac{g^2}{F^2} \, dt }
-1
\right]
\frac{g(t_1) }{F(t_1)}
\left[\alpha
-
\int_0^{t_1}
\frac{ g} {F} \,dt
\right]
\]
\[
+
\frac{
\int_0^{t_1 }
H\,
\frac{g^2}{F^2} \, dt
}
{\int_{0}^{t_1} H  \frac{g^2}{F^2}  \, dt
- H(t_1) \, \int_0^{t_1 }
\frac{g^2}{F^2} \, dt }
\]
\[
=
\frac{
\frac{g(t_1) \, H(t_1) }{F(t_1)}
\left[\alpha
-
\int_0^{t_1}
\frac{ g} {F} \,dt
\right]
+
\int_0^{t_1 }
H\,
\frac{g^2}{F^2} \, dt
}
{\int_{0}^{t_1} H  \frac{g^2}{F^2}  \, dt
- H(t_1) \, \int_0^{t_1 }
\frac{g^2}{F^2} \, dt }
\]
%
%
Hence we have
\[ \lambda_1 \, \frac{g(t_1)}{F(t_1)}  = 1 + \gamma(t_1) \, \frac{g(t_1) \, H(t_1)}{F(t_1)}.\]

Assume that $g$ is continuously differentiable and we have
\begin{equation}
\label{derivativeassumption}
g'(t) \leq f(t) \, g(t).
\end{equation}
Assumption
(\ref{derivativeassumption}) implies that the function
$\frac{g}{F}$ is decreasing.

Assumption
(\ref{derivativeassumption}) implies that the function
$\frac{g}{F}$ is decreasing.
Since the function $(\lambda_1 - \gamma(t_1) \, H)$ is decreasing
and $\frac{g}{F}>0$  this implies that also the product
\[\frac{g}{F} \, (\lambda_1 - \gamma(t_1) \, H)\]
is decreasing as a function of time.

%
Then the  optimal control
$\hat u$ as defined in (\ref{optimalcontroldefinition})
(with $\lambda = \lambda_1$ and $\gamma = \gamma(t_1)$)
is decreasing,
$\hat u(t_1)=0$
and the
support of the optimal control
$\hat u$   
is contained in $[0, t_1]$.
With $\lambda_1$ defined as in  (\ref{lambdadefinition}),
equation (\ref{momentequationt0}) holds.
%
This implies that the optimal control
$\hat u$ as defined in (\ref{optimalcontroldefinition})
satisfies  (\ref{lambdawahlbedingung}).
Thus we have shown the following statement:

{\em
If (\ref{derivativeassumption}) holds,
for all $t_0 \in (0, T)$
such that
(\ref{t1condition}) holds (with $t_1 = t_0$)
there is
a weight $\gamma>0$ such that the support of
the corresponding optimal control is
contained in  $[0, \, t_0]$.
}

Note that in Example \ref{bsp1},
we have $f(t)=1$
and $g(t) = \exp(t) = g'(t)$,
hence (\ref{derivativeassumption}) holds.
This explains why  in the first example,
for sufficiently large values of $T$ no
terminal constraint is necessary.
As a second example, for
the constant function $g(t)=1$,
(\ref{derivativeassumption}) also holds.

\section{General results in Hilbert spaces}

In this section, we study
optimal control  problems in a Hilbert space setting.
In this way, we obtain results that
we can apply to systems that are governed by
partial differential equations.
Let $X$ and $U$ be Hilbert spaces with the inner products
$\langle \cdot, \cdot\rangle_X$,
$\langle \cdot, \cdot\rangle_U$
and the corresponding norms $\|\cdot\|_X$,  $\|\cdot\|_U$ respectively.
We use $T>0$ to denote the terminal time
of our optimal control problems.
The space $X$ contains the current state and
the space $U$ is used as a framework for the
control functions in $L^2(0, T;U)$.

Let $A: {\cal D}(A)\subset X \rightarrow X$
be the generator of a strongly continuous
semigroup, and let $B$ denote an admissible control operator.
As in \cite{tu:obser}, Proposition 4.2.5., we consider control systems
of the form
\begin{equation}
\label{1}
\left\{
\begin{array}{rcl}
x' + A x & = & B u,
\\
x(0) & = & x_0
\end{array}
\right.
\end{equation}
where $x_0 \in X$ is a given initial state.
For all $u\in L^2(0, T;U)$,
the Cauchy problem (\ref{1})
has a unique solution $x \in C([0, T]; X)$
(see \cite{phillips}).



\subsection{Exact controllability}  
Assume
that (\ref{1}) is
{\em exactly controllable
using $L^2$--controls} in time $t_0  > 0$, that is there
exists a constant $C_1>0$ such that for all initial states $x_0\in
X$
and all terminal  states $x_1\in X$
there is a control $u\in L^2(0,t_0;U)$ such that the solution $x\in
C([0,t_0];X)$ of (\ref{1})
satisfies
\begin{equation}
\label{2}
\left\{
\begin{array}{rcl}
x(t_0)  & = & x_1.
\\
\|u\|_{L^2(0,t_0 ;U)} & \leq  & C_1 (\| x_0\|_X + \|x_1\|_X).
\end{array}
\right.
\end{equation}


%
%


Let a desired state $x_d \in X$ be given.
Due to the exact controllability assumption,
there exists a control $u_{\rm exact} \in L^2(0, t_0; \, U)$
such that the solution $x_{\rm exact} \in
C([0,t_0];X)$ of (\ref{1})
satisfies
\begin{equation}
\begin{array}{rcl}
x_{{\rm exact}}(t_0)  & = & x_d.
\end{array}
\end{equation}
We assume that $x_d$ is a {\em holdable state} in
the sense that
we can extend $u_{\rm exact}$ to the time interval
$[0, \, T]$  by a constant control
$u_d$
on $[t_0,\, T]$ such that
for the corresponding state
for all $t\in (t_0,\, T)$
we have the equation
$x_{{\rm exact}}(t)  =  x_d$
and
$u_{{\rm exact}}(t)  =  u_d$.
Thus on the time--interval   $(t_0,\, T)$
we have $A\, x_{\rm exact} = A \, x_d = B\, u_{\rm exact}$.

\subsection{An optimal control problem with $\max$-norm penalization}

First we consider a tracking term
with the maximum-norm.
For systems that are exactly controllable,
the optimal control steers the system
 to the desired state after the prescribed time $t_0$.

For
$\gamma>0$
we consider the following optimization problem:
\[
{\bf P}(T, \, \gamma)
\left\{
\begin{array}{lrr}
\min\limits_{u\in L^2(0,t; \, U)}
\frac{1}{2} \; \|u - u_d \|^2_{L^2(0,T;U)}
+
\gamma \,
\,
\max_{t \in [t_0,\, T]}
\| x(s) - x_d \|_X
\\
\mbox{\rm subject to} \\
 x' + A x = Bu, \; x(0)=x_0.
\end{array}
\right.
\]
In problem ${\bf P}(T, \, \gamma)$
 the end condition $x(t)=x_d$ does not appear.
Note that problem ${\bf P}(T,\gamma)$ has a unique solution.

Our goal is to show that, due to the property of
{exact controllability using $L^2$--controls}
of the system, for $\gamma$ sufficiently large
the optimal  state $x_T$ satisfies the  condition
\[
x_T(t)=x_d
\]
for all $t\in [t_0,\, T]$.
%
%
%
A precise statement is given in the following theorem:
\begin{theorem}
\label{allgemein}
Assume that
$T > t_0$
 and that
the system (\ref{1})
is exactly controllable.
 If
$
\gamma >  0
$
is sufficiently large, for all $t\in [t_0,\, T]$
 the solution
 $(u_T,\, x_T)$
 of problem ${\bf P}(T,\, \gamma)$
satisfies the equation
 \[x_T(t) = x_d.\]
\end{theorem}
{\bf Proof:}
An application of the Direct Method of the Calculus of Variations
shows that
a solution
of  ${\bf P}(T, \gamma)$  exists.
The strict convexity of the control cost
$\frac{1}{2}\,\| \cdot\|^2_{L^2(0,T;U)}$
implies that the solution of  ${\bf P}(T, \, \gamma)$   is uniquely determined.
%
Choose
\begin{equation}
\label{contra12019}
\gamma > C_1 \,  \|u_{\rm exact} - u_d\|_{L^2(0,t_0;U)}.
\end{equation}

Similarly as in \cite{guzuexact}, consider the optimal control  problem
\[
{\bf Q}(T,\gamma)
\left\{
\begin{array}{lrr}
\min\limits_{u\in L^2(0,T; \, U)}
\frac{1}{2} \; \|u - u_d \|^2_{L^2(0,T;U)}
+
\gamma \,
\,
\| x(t_0) - x_d \|_X
\\
\mbox{\rm subject to} \\
 x' + A x = Bu, \; x(0)=x_0.
\end{array}
\right.
\]
Let $(u^\ast, \, x^\ast)$ denote the solution of
${\bf Q}(T,\gamma)$.
Now similarly as in Theorem 1 in (\cite{guzuexact}),
we show that
$ x^\ast(t_0)= x_d$ by an indirect proof.

Suppose that $ x^\ast(t_0) \not= x_d$.
Then the objective functional of ${\bf Q}(T,\gamma)$ is differentiable at
$(u^\ast, \, x^\ast)$
and  the necessary optimality conditions imply
\begin{equation}
\label{ten19}
\int_{0}^{T}
 \langle u^\ast - u_d, v \rangle_U  \, dt
+
\gamma \,
\frac{ \langle x^\ast(t_0) - x_d, y \rangle_X}
 {\|x^\ast(t_0) - x_d\|_X}
=0
\end{equation}
for all $v\in L^2(0,t_1;U)$ where $y$ solves
\[
y'+ Ay =
Bv, \; y(0)=0.
\]

Due to the exact controllability of the system,
we can choose a control $\tilde v\in L^2(0, t_0;U)$
such that for the corresponding state $\tilde y$ we have
\[\tilde y(t_0) = \frac{ x^\ast(t_0) - x_d }{\|x^\ast(t_0) - x_d\|_X}\]
and
\begin{equation}
\label{vschlangeungleichung}
\|\tilde v\|_{L^2(0, t_0; U)} \leq C_1.
\end{equation}
We extend $\tilde v$ to an element of $L^2(0, T;U)$
by the definition $\tilde v(s) = 0$ for all $ s\in (t_0,\, T)$.
Then the necessary optimality condition yields the equation
\begin{equation}
\label{ten19a}
\int_{0}^{T}
 \langle u^\ast - u_d, \tilde v \rangle_U  \, dt
+
\gamma \,
\frac{ \langle x^\ast(t_0) - x_d,
\frac{ x^\ast(t_0) - x_d }{\|x^\ast(t_0) - x_d\|_X}
 \rangle_X}
 {\|x^\ast(t_0) - x_d\|_X}
=0.
\end{equation}

This implies
the equation
\begin{equation}
\label{gleichung27092019}
\left|
\int_{0}^{T}
 \langle u^\ast - u_d, \tilde v \rangle_U  \, dt
\right|
=
\gamma.
\end{equation}
On the other hand, we have the inequality
\[
\left|
\int_{0}^{T}
 \langle u^\ast - u_d, \tilde v \rangle_U  \, dt
\right|
\leq \|u^\ast - u_d\|_{L^2(0,T; U)} \,  \|\tilde v\|_{L^2(0,T; U)}.
\]
Since the control
$u_{\rm exact}$ is feasible for ${\bf Q}(T,\gamma)$,
we have the inequality
\[ \frac{1}{2} \, \|u^\ast - u_d\|_{L^2(0,T; U)}^2
\leq  \frac{1}{2} \,  \|u_{\rm exact}  - u_d\|_{L^2(0,T; U)}^2 +
\gamma \, \| y_{\rm exact}(t_0)  - y_d\|_X
\]
\[
=\frac{1}{2} \,  \|u_{\rm exact}  - u_d\|_{L^2(0,T; U)}^2.
\]
Hence
\[  \|u^\ast - u_d\|_{L^2(0,T; U)}
\leq
\|u_{\rm exact}  - u_d\|_{L^2(0,T; U)}.
\]
Moreover, (\ref{vschlangeungleichung}) implies
\[\|\tilde v\|_{L^2(0, T; U)} \leq C_1. \]
Hence
(\ref{gleichung27092019})
implies
\[\gamma \leq C_1 \, \|u_{\rm exact}  - u_d\|_{L^2(0,T; U)},\]
which is a contradiction to
(\ref{contra12019}).
Thus we have shown that
$ x^\ast(t_0)= x_d$.
This implies that for
$s\in (t_0,\, T]$ we have
$u^\ast(s) = u_d$ and $x^\ast(s)= x_d$.

Let $v_Q$ denote the optimal value of
${\bf Q}(T,\gamma)$ and $v_P$
denote the optimal value of
${\bf P}(T,\gamma)$.
Then the definition of the corresponding objective  functionals
implies the inequality
\[v_Q \leq v_P.\]

Since the control $u^\ast$ is feasible for ${\bf P}(T, \gamma)$,
we also have the inequality
\[
v_P =
\frac{1}{2} \; \|u_T - u_d \|^2_{L^2(0,T;U)}
+
\gamma \,
\,
\max_{t \in [t_0,\, T]}
\| x_T(s) - x_d \|_X
\]
\[
\leq
\frac{1}{2} \; \|u^\ast - u_d \|^2_{L^2(0,T;U)}
+
\gamma \,
\,
\max_{t \in [t_0,\, T]}
\| x^\ast(s) - x_d \|_X
\]
\[
=
\frac{1}{2} \; \|u^\ast - u_d \|^2_{L^2(0,t_0;U)}
= v_Q.
\]
Thus we have $v_P = v_Q$,
and $(u^\ast, \, x^\ast)$ is
an optimal control/state pair for ${\bf Q}(T,\gamma)$.
Since the solution is unique, this implies
the assertion.
$\Box$.


\subsection{An optimal control problem for nodal profile exactly controllable systems}
\label{nodal}

Motivated by application problems in the operation of gas piplines,
the exact controllability of nodal profiles has
been introduced in \cite{controllabilitynodalprofile}, see also
\cite{libook}.
The assumption of  exact controllability of nodal profiles
also allows to derive a result about the exactness of an $L^2$-norm penalty term.

Let a Hilbert space $Z$, $t_0 \in (0, \, T)$  and
 a linear map
$\Pi :L^2(0, T;X) \rightarrow L^2(t_0, T;Z)$ be given.
In the applications, typically $\Pi$ will
be some trace operator, for example
the boundary trace of the system state
restricted to the time-interval $[t_0,\, T]$,
see \cite{controllabilitynodalprofile}.

Assume
that (\ref{1}) is
{\em nodal profile exactly controllable
using $L^2$--controls}
in time $t_0  > 0$, that is there
exists a constant $C_1>0$ such that for all initial states $x_0\in
X$
and all nodal profiles
$z \in L^2(t_0, T; Z)$
there is a control $u\in L^2(0, T;U)$ such that the solution $x\in
C([0, T];X)$ of (\ref{1})
satisfies for all $t\in [t_0,\, T]$
\begin{equation}
\label{2nodal}
\left\{
\begin{array}{rcl}
\Pi x(t)  & = & z(t),
\\
\|u\|_{L^2(0, T ;U)} & \leq  & C_1  \, ( \| x_0\|_X + \|z\|_{L^2(t_0, T; Z)}).
\end{array}
\right.
\end{equation}

\begin{remark}
The exact boundary controllability of nodal profile for hyperbolic systems
is discussed in \cite{libook}.
\end{remark}

For
$\gamma>0$
we consider the following optimization problem:
\[
{\bf S}(T, \, \gamma)
\left\{
\begin{array}{lrr}
\min\limits_{u\in L^2(0,t; \, U)}
\frac{1}{2} \; \|u - u_d \|^2_{L^2(0,T;U)}
+
\gamma \,
\,
\sqrt{
\int\limits_{t_0}^T
\| \Pi x(s) -  \Pi x_d \|_Z^2
\, ds
}
\\
\mbox{\rm subject to} \\
 x' + A x = Bu, \; x(0)=x_0
\end{array}
\right.
\]
where as before, $x_d \in X$ is the desired holdable state.
In problem ${\bf S}(T, \, \gamma)$
 the end condition $x(t)=x_d$ does not appear.
Note that problem ${\bf S}(T,\gamma)$ has a unique solution.

\begin{remark}
Optimization problems of a similar structure
with a differentiable tracking term have been considered in
\cite{guha} and \cite{gupafa}.
\end{remark}

Due to the nodal profile  exact controllability assumption,
there exists a control $v_{\rm exact} \in L^2(0, t_0; \, U)$
such that the solution $p_{\rm exact} \in
C([0,t_0];X)$ of (\ref{1})
satisfies
\begin{equation}
\begin{array}{rcl}
\Pi p_{{\rm exact}}(t)  & = & \Pi x_d
\end{array}
\end{equation}
for all $t\in [t_0, \, T]$.

Our goal is to show that, due to the property of
{nodal profile exact controllability using $L^2$--controls}
of the system, for $\gamma$ sufficiently large
the optimal  state $x_T$ satisfies the  condition
\[
\Pi x_T(t)= \Pi x_d
\]
for all $t\in [t_0,\, T]$.
In the application  in supply systems,
this means
that on the time interval $[t_0,\, T]$,
the nodal profile that is
desired by the customer is attained exactly.
%
%
A precise statement is given in the following theorem:
\begin{theorem}
\label{allgemeinnodal}
Assume that
$T > t_0$
 and that
the system (\ref{1}) is
nodal profile exactly controllable.
 If
$
\gamma >
 C_1 \,  \|u_{\rm exact} - u_d\|_{L^2(0,t_0;U)}$,
for all $s\in [t_0,\, T]$  the solution
 $(u_T,\, x_T)$
 of problem ${\bf S}(T,\, \gamma)$
satisfies the equation
 \[\Pi x_T(s) =  \Pi  x_d.\]
\end{theorem}
{\bf Proof:}
An application of the Direct Method of the Calculus of Variations
shows that
a solution
of  ${\bf S}(T, \gamma)$  exists.
The strict convexity of the control cost
$\frac{1}{2}\,\| \cdot\|^2_{L^2(0,T;U)}$
implies that the solution of  ${\bf S}(T, \, \gamma)$   is uniquely determined.
%
Choose
\begin{equation}
\label{contra12019nodal}
\gamma > C_1 \,  \|v_{\rm exact} - u_d\|_{L^2(0,t_0;U)}.
\end{equation}

Suppose that there exists $\tau \in [t_0,\, T]$ such that
$ \Pi x^\ast(\tau) \not =  \Pi x_d$.
Then $ \|\Pi x^\ast  -  \Pi x_d\|_{L^2(t_0, T; Z)} \not = 0$.
Hence the objective functional of ${\bf S}(T,\gamma)$ is differentiable
in $(u^\ast, \, x^\ast)$
%
and  the necessary optimality conditions imply
\begin{equation}
\label{ten190510}
\int_{0}^{T}
 \langle u^\ast - u_d, v \rangle_U  \, dt
+
\gamma \,
\int_{t_0}^{T}
\frac{ \langle \Pi x^\ast(\tau ) - \Pi x_d,  \Pi y(\tau) \rangle_Z}
 { \|\Pi x^\ast   - \Pi x_d\|_{L^2(t_0, T; Z)} }
 \, d\tau
=0
\end{equation}
for all $v\in L^2(0,T;U)$ where $y$ solves
\[
y'+ Ay =
Bv, \; y(0)=0.
\]

Due to the nodal profile exact controllability of the system,
we can choose a control $\tilde v\in L^2(0, t_0;U)$
such that for the corresponding state $\tilde y$ we have
for all $\tau \in [t_0, \, T]$
\[\Pi \tilde y(\tau) = \frac{ \Pi x^\ast(\tau) -  \Pi x_d }{ \|\Pi x^\ast   - \Pi x_d\|_{L^2(t_0, T; Z)} }\]
and
\[\|\tilde v\|_{L^2(0, T; U)} \leq C_1. \]

Then the necessary optimality condition
(\ref{ten190510})
 yields the equation
\begin{equation}
\label{ten19a2019}
\int_{0}^{T}
 \langle u^\ast - u_d, \tilde v \rangle_U  \, dt
+
\gamma \,
\int_{t_0}^{T}
\frac{ \langle \Pi x^\ast(\tau ) - \Pi x_d,
\frac{ \Pi x^\ast(\tau ) -  \Pi x_d }{ \|\Pi x^\ast   - \Pi x_d\|_{L^2(t_0, T; Z)}  }
 \rangle_Z}
 { \|\Pi x^\ast   - \Pi x_d\|_{L^2(t_0, T; Z)}  }
 \, d\tau
=0.
\end{equation}

This implies
the equation
\begin{equation}
\label{gleichung27092019nodal}
\left|
\int_{0}^{T}
 \langle u^\ast - u_d, \tilde v \rangle_X  \, dt
\right|
=
\gamma.
\end{equation}
On the other hand, we have the inequality
\[
\left|
\int_{0}^{T}
 \langle u^\ast - u_d, \tilde v \rangle_U  \, dt
\right|
\leq \|u^\ast - u_d\|_{L^2(0,T; U)} \,  \|\tilde v\|_{L^2(0,T; U)}.
\]
Since the control
$v_{\rm exact}$ is feasible for ${\bf S}(T,\gamma)$,
we have the inequality
\[ \frac{1}{2} \, \|u^\ast - u_d\|_{L^2(0,T; U)}^2
\leq  \frac{1}{2} \,  \|v_{\rm exact}  - u_d\|_{L^2(0,T; U)}^2 +
\gamma \,\int_{t_0}^T \| \Pi p_{\rm exact}(\tau )  - \Pi x_d\|_Z \, d\tau
\]
\[
=\frac{1}{2} \,  \|v_{\rm exact}  - u_d\|_{L^2(0,T; U)}^2.
\]
Hence
\[  \|u^\ast - u_d\|_{L^2(0,T; U)}
\leq
\|v_{\rm exact}  - u_d\|_{L^2(0,T; U)}.
\]
Moreover, we have
\[\|\tilde v\|_{L^2(0, T; U)} \leq C_1. \]
Hence
(\ref{gleichung27092019nodal})
implies
\[\gamma \leq C_1 \, \|v_{\rm exact}  - u_d\|_{L^2(0,T; U)},\]
which is a contradiction to
(\ref{contra12019nodal}).
Thus we have shown that
$ \Pi x^\ast  = \Pi x_d$
on $[t_0,\, T]$.
%
%
This implies
the assertion.
$\Box$.


\subsection{An optimal control problem with $L^1$-norm tracking term}
In this section we present a result
about the
finite-time turnpike structure of the optimal state and the optimal control
that we have shown
under the assumption of exact controllability (\ref{2})
for an optimal control problem with an
 $L^1$-norm tracking term with a singular weight
 in the objective function.


For
$\gamma>0$
we consider the following optimal control problem
${\bf R}(T ,\gamma)$
with $L^1$-norm tracking term:
\[
{\bf R}(T ,\gamma)
\left\{
\begin{array}{lrr}
\min\limits_{u\in L^2(0, T ; \, U)}
\frac{1}{2} \; \|u - u_d \|^2_{L^2(0, T;U)} +
\gamma
\,
 \int_{t_0}^T
 \frac{1}{s - t_0}
\| x(s) - x_d \|_X \, ds
\\
\mbox{\rm subject to} \\
 x' + A x = Bu, \; x(0)=x_0.
\end{array}
\right.
\]
In problem ${\bf R}(T  ,\gamma)$
 the end condition $x(T  )=x_d$ does not appear.
Problem ${\bf R}(T ,\gamma)$ has a unique solution.

Our goal is to show that, due to the property of
{exact controllability using $L^2$--controls}
of the system, for $\gamma$ and $T$ sufficiently large
the optimal  state $x_T$ for ${\bf R}( T  ,\gamma)$ satisfies the  condition
\[
x_T(t)=x_d
\]
for all $t\in (t_0, \, T]$.
%
A precise statement is given in the following theorem:
\begin{theorem}
\label{allgemeinI}
Assume that
$T > t_0$
 and that
the system (\ref{1})
is exactly controllable.
 If
$
\gamma >  0,
$
 the solution
 $(u_T,\, x_T)$
 of problem ${\bf R}(T,\, \gamma)$
satisfies the equation
 \[x_T(t) = x_d\]
 for all $t\in [t_0,\, T]$.
\end{theorem}
{\bf Proof:}
Since $u_{\rm exact}$ is a  feasible control for ${\bf R}(T, \gamma)$,
evaluating the objective function of  ${\bf R}(T, \gamma)$ at $ u_{\rm exact}$
yields the inequality
\begin{equation}
\label{upper22sequentiell}
\|u_{T} - u_d\|^2_{L^2(0,{T};U)}
\leq
\|u_{T} - u_d\|^2_{L^2(0,t_0;U)} + 2\,
 \gamma
 \,\int_{t_0}^{T}
 \frac{1}{s- t_0}
 \|x_{{\rm exact}}(s) - x_d\|_X \, ds
 \end{equation}
 \[
 =  \|u_{\rm exact} - u_d\|^2_{L^2(0,t_0;U)} .
\]

An application of the Direct Method of the Calculus of Variations
shows that
a solution
of  ${\bf R}(T, \gamma)$  exists.
For the optimal control/state pair we use the notation
$(u_T,\, x_T)$.

 If there exists $\hat t\in (0, T)$ with
 $x_T( \hat t)=x_d$,
 the optimal way to continue
 the control for $s\in (\hat t, \, T]$ is with
 $(u_d, \, x_d)$,
 hence for all $s\in (\hat t,\, T]$ we have
 $x_T( s)=x_d$.
%
%
%
%
%
%
%

Suppose that there exists a number
$t_1 \in (t_0, \, T]$ such that
$x_{T}(t_1)\not=x_d$.
Then for all $t\in [t_0,\, t_1)$, we also have
 $x_{T}(t)\not=x_d$.
In particular,   for all $t\in [t_0,\, t_1]$, we have
 $\|x_{T}(t) - x_d\|_X >0$.
 Since $x_T$ is continuous, this implies that
 \[\inf_{t \in  [t_0,\, t_1]} \|x_{T}(t) - x_d\|_X = \varepsilon >0.\]
 This implies
 \[
  \int_{t_0}^{t_1}
 \frac{1}{s - t_0}
\| x(s) - x_d \| \, ds
\geq \varepsilon
 \int_{t_0}^{t_1}
 \frac{1}{s - t_0}
=
\infty.
\]
Hence $x_T$ cannot be optimal, and this is a contradiction.
$\Box$

\section{Examples}
\label{example}
In this section we present some examples to
illustrate our results about the finite-time turnpike
phenomenon.
We start with one example
with a system that is governed by an ordinary
differential equation and then we present
examples with partial differential equations.

\begin{example}
\label{beispielnumerik}
Let us first return to Example \ref{bsp1}.
Here we present numerical results
that illustrate that the numerical solution
for the discretized optimal control
problem
where for $T=2$ the interval $[0,2]$
has been replaced with  a grid of
201 equidistant points and
 the ordinary differential equation  has been replaced by
a discrete time-system with the Euler backwards
discretization.

The resulting  optimization problem
has been solved numerically with a
standard method from matlab.
To improve the performance, in
the numerical experiments the constraints
$u\geq 0$ and $y\leq 0$ have been included in
the problem.
(As shown in Example \ref{bsp1}, they do not change the solution).
The numerical results
are presented in Figure \ref{fig_solution1} for $\gamma = \tfrac{1}{2}$,
Figure \ref{fig_solution2}  for $\gamma =1$ and Figure \ref{fig_solution3} for $\gamma =2$.

\begin{figure}[hbt]
\includegraphics[width=0.95\textwidth]{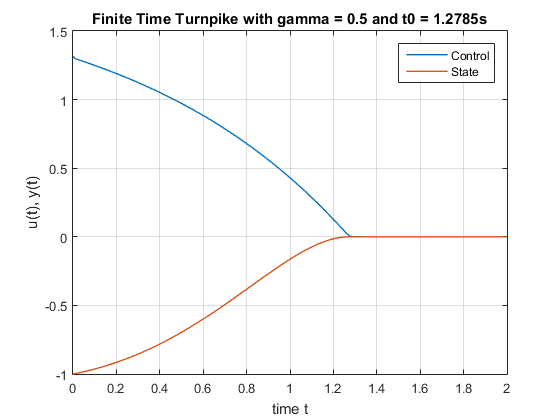}
\caption{
\label{fig_solution1}
The figure shows the optimal control and the optimal state
 as approximate solutions of problem ${\bf (OC)}_T$ for  $T=2$ and $\gamma=\tfrac{1}{2}$
 defined in Example \ref{beispielnumerik}.}
\end{figure}

\begin{figure}[hbt]
\includegraphics[width=0.95\textwidth]{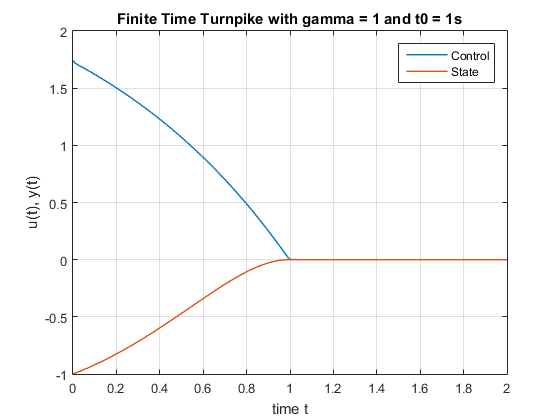}
\caption{
\label{fig_solution2}
The figure shows the optimal control and the optimal state
 as approximate solutions of problem ${\bf (OC)}_T$ for  $T=2$ and $\gamma=1$
  defined in Example \ref{beispielnumerik}.}
\end{figure}

\begin{figure}[hbt]
\includegraphics[width=0.95\textwidth]{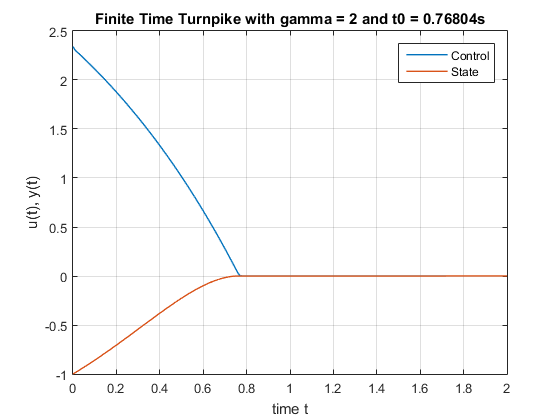}
\caption{
\label{fig_solution3}
The figure shows the optimal control and the optimal state
 as approximate solutions of problem ${\bf (OC)}_T$ for  $T=2$ and $\gamma=2$
  defined in Example \ref{beispielnumerik}.}
\end{figure}

\end{example}


Now we present examples of optimal control problems
where  Theorem  \ref{allgemein}
or   Theorem  \ref{allgemeinI}  is applicable.
These theorems assume that the system is exactly controllable.

\begin{example}
Now we consider a problem of  optimal torque  control for an Euler--Bernoulli beam.
%
Let $y_0 \in H^2(0, 1)$ and $y_1\in H^1(0, 1)$ be given.
We study the following optimal control problem:
\[
\left\{
\begin{array}{l}
\min\limits_{u\in L^2(0,T)} \;
\frac{1}{2} \|u^2(t)\|^2\, dt
+
\gamma \,
\max_{ t \in [t_0,\,T] }
\| y(t,\, \cdot )\|_{L^2(0,1)}
\; \,{\sf subject } \;{\sf to}\;
\\
y(0,x)= y_0(x),\; y_t(0,x)= y_1(x),\; x\in (0,1)
\\
y(t,0)= 0,\; {y_{xx}(t,0)= u(t),}\; t\in (0,T)
\\
y(t,\, 1) = y_{xx}(t,1) =0,
\\
y_{tt}(t,x)= -  y_{xxxx}(t,x),\;(t,x)\in (0,T) \times (0,1).
\end{array}
\right.
\]
We have $U= L^2(0,1)$ and
$X=L^2(0,1)$.
Note that the Euler--Bernoulli beam is exactly
controllable in arbitrarily short times
(see \cite{tu:obser}, Example 11.2.8),
so in this case
$ t_{0}>0$ can be chosen arbitrarily small.
 Theorem  \ref{allgemein} implies
 that if $\gamma$ is chosen sufficiently large
 the beam is steered to a position of rest in
 the time $t_0>0$,
\end{example}

\begin{example}
Consider the  problem of
optimal Neumann boundary
control of the wave equation.
Define $Q=(0, \, T) \times (0,1)$.
Here we have
  $U = L^2(0,\, 1)$,
  $X=L^1(0,1)$,

Let $y_d \in X$ and  $ u_d \in U$  be given.
Consider the optimal control problem
\[
\left\{
\begin{array}{l}
\min\limits_{u\in U } \;
\frac{1}{2}
\int\limits_0^T \, \left( u(t) - u_d \right)^2  \, dt
+
\int\limits_2^T \,\frac{1}{t-2}\,
\int\limits_0^1
\left| y(t,\, x) - y_d \right| \,dx\,dt
\; \,{\sf subject } \;{\sf to}\;
\\
y(0,x)= 0,\; y_t(0,x)= 0,\; x\in (0,1)
\\
y(t,0)= 0,\; {y_x(t,1)= u(t),}\; t\in (0,T)
\\
y_{tt}(t,x) -  y_{xx}(t,x) = 0, \;
(t,x)\in Q.
\end{array}
\right.
\]

Our results show that
the solution
has a turnpike structure
as described in Theorem  \ref{allgemeinI}.
The optimal control problem is similar to the
Neumann optimal boundary control problem
with a differentiable objective function considered
in
\cite{guzu}.

\end{example}

Now we present an example where
Theorem \ref{allgemeinnodal} is applicable,
that assumes that the system is nodal profile exactly controllable.

\begin{example}
Now we consider a problem or  optimal control
where Theorem \ref{allgemeinnodal} is applicable.
The problem is similar as in \cite{guha}, but
in the tracking term instead of the squared $L^2$-norm
we take the $L^2$-norm.
The motivation for this type of problem
where the boundary trace of the state is driven to a desired profile
comes from the operation of  networks of gas pipelines,
where the aim is to satisfy customer demands
in an optimal way.

We consider a $2\times 2 $ system in diagonal form.
Let a length $L>0$ and a time interval $[0,T]$ be given.
Let  $d_-$ and $d_+$ be real numbers
such that
\[d_- < 0 < d_+.\]
Define the
diagonal matrices
\[
D =
\left(
\begin{array}{cc}
d_+&  0
\\
0   &  d_-
\end{array}
\right).
\]
For all $x\in [0,\, L]$, let $M(x)$ denote a
$2\times 2$ matrix
that depends continuously on $x$.
Assume that for all $x\in [0,\, L]$ the matrix
$M(x)$ is  positive semi--definite.
Let $\eta_0\leq 0$ be  a real number.
%

Consider the linear hyperbolic partial differential equation
\begin{equation}
\label{pde}
r_t + D\, r_x =
\eta_0\, M\, r
\end{equation}
where for $x\in (0,\, L)$ and $t\in (0,T)$, the state is given by
$r(t,\,x)= \left(
\begin{array}{r}
r_+(t,\,x)
\\
r_-(t,\,x)
\end{array}
\right).
$

Let real numbers
$R_+^{d}$ and $R_-^{d}$ be given.
To obtain an initial boundary value
problem, in addition to (\ref{pde})
we consider the initial condition
$r(0,\,x)=0$ for $x\in (0,\, L)$ at the time $t=0$
and
for $t\in (0,T)$
the Dirichlet
 boundary conditions
$
r_+(t,\,0)= u_+(t),\;
r_-(t,\, L)=  R^d_-,
$
with a boundary control
$u_+$  in  $L^2(0,T)$.
%
%
The resulting initial boundary value problem
\begin{equation}
\label{linearizedsystem}
\left\{
\begin{array}{l}
r(0,\,x)=0,
\\
r_t + D\, r_x = \eta_0\,M\, r,
\\
r_+(t,\,0)= u_+(t),
\\
r_-(t,\, L)= R^d_-
\end{array}
\right.
\end{equation}
has a solution
$r\in C([0,T], L^2((0,\, L); {\mathbb R}^2))
$.
Moreover, for the boundary traces
of the solution we have
$r_+(\cdot, L)$, $r_-(\cdot,\,0)\in L^2(0,\, T)$.

For $x= (x_+,\, x_-) \in {\mathbb R}^2$,
we use the notation $\| x \|_{{\mathbb R}^2} = \sqrt{x_+^2 + x_-^2}$.
For $u=(u_+,u_-) \in ( L^2(0,\, T))^2$ and $R=(R_+,R_-)\in ( L^2(0,\, T))^2$,
define the objective function
\[J(u,\, R)\]
\begin{equation}
\label{objectiveneu}
=\int_0^T
\tfrac{1}{2} \,
(u_+(t))^2 \, dt
+ \gamma \,
\int_{t_0}^{T- t_0}
\| (R_+(t)- R_+^{d} ,\,  R_-(t)- R_-^{d}  \|_{{\mathbb R}^2} \, dt
.
\end{equation}
Then if $L$ is sufficiently small
and $T$
and $t_0 < T$ are sufficiently large, the system
is nodal profile exactly controllable
and  Theorem \ref{allgemeinnodal} is applicable for the optimal control problem
\begin{equation}
\label{ocplambda}
\left\{
\begin{array}{l}
\min_{u_+
\in  L^2(0,\, T)} J(u_+,\, (r_+(\cdot, L),\,r_-(\cdot,\,L))  )\;
\\
\mbox{\rm subject to
(\ref{linearizedsystem})}.
\end{array}
\right.
\end{equation}
In fact the result of Theorem \ref{allgemeinnodal}
can be interpreted as a {\em finite-time turnpike result}
(or exact turnpike), where
the system is driven to a desired stationary state in finite time

\end{example}

\section{Conclusion}
\label{conclusion}
We have shown that a finite-time turnpike phenomenon
occurs for  problems of optimal control
with nondifferentiable norm tracking terms.

We have first considered systems
that are governed by ordinary differential equations.
In the objective functions, $L^1$-norm tracking terms
are used.
The finite-time turnpike means that
after finite time
the optimal state reaches the
desired state.
For infinite-dimensional systems,
we have shown that a finite-time turnpike phenomenon
occurs for  problems of optimal control
for systems that are exactly controllable
with a $\max$-norm type  tracking term
and a weighted $L^1$-norm tracking term.
For systems that  are nodal profile exactly controllable,
we have shown that
a finite-time turnpike phenomenon
occurs with an $L^2$-norm tracking term.
\\
This work was supported by the DFG grant CRC/Transregio 154, project C03
and C05.

\end{document}